\documentclass[12pt,a4paper]{amsart}
\usepackage[utf8]{inputenc}
\usepackage{amsfonts}
\usepackage{amssymb}
\usepackage{amsmath}
\usepackage{amsthm}
\usepackage{cite}
\usepackage{enumitem}
\usepackage{tikz}
\usepackage{subfigure}
\usepackage{dsfont}

\theoremstyle{plain}
\newtheorem{thm}{Theorem}

\newtheorem{lem}{Lemma}
\newtheorem{prop}{Proposition}

\newtheorem{rmk}{Remark}

\providecommand{\N}{\mathbb{N}}
\providecommand{\R}{\mathbb{R}}
\providecommand{\Z}{\mathbb{Z}}

\newcommand{\norm}[1]{\left\lVert #1 \right\rVert}
\newcommand{\skp}[2]{\left< #1 ,#2 \right>}

\renewcommand{\i}{\mathrm{i}}

\makeatletter
\newcommand{\pushright}[1]{\ifmeasuring@#1\else\omit\hfill$\displaystyle#1$\fi\ignorespaces}
\makeatother

\renewcommand{\b}[1]{\mathbf{#1}}

\newenvironment{proofof}[1]{
\vspace*{0.4cm}
\textbf{Proof of \textsc{#1}.} \\
}{\qed}

\newcounter{stepcount}
\newcounter{substepcount}[stepcount]

\newenvironment{steps}[0]{
\setcounter{stepcount}{0}
\setcounter{substepcount}{0}
\newcommand{\step}[1]{\stepcounter{stepcount}
\vspace*{0.2cm} \underline{\textsc{Step} \arabic{stepcount}}: \textit{##1} 

\vspace*{0.1cm}
}

}{ }

\setitemize{itemsep=+2pt}
\setenumerate{itemsep=+2pt}
\begin{document}

\allowdisplaybreaks

\title{Breather Solutions of the Cubic Klein-Gordon Equation}

\author{Dominic Scheider}
\address{ D. Scheider\hfill\break
Karlsruhe Institute of Technology \hfill\break
Institute for Analysis \hfill\break
Englerstra{\ss}e 2 \hfill\break
D-76131 Karlsruhe, Germany}
\email{dominic.scheider@kit.edu}
\date{\today}


\begin{abstract}
We obtain real-valued, time-periodic and radially symmetric solutions of the cubic Klein-Gordon equation
$$
	\partial_t^2 U - \Delta U + m^2 U = \Gamma (x) U^3
	\quad \text{on } \R \times \R^3,
$$
which are weakly localized in space.
Various families of such ``breather'' solutions are shown to bifurcate from any given nontrivial stationary solution. The construction of weakly localized breathers in three space dimensions is, to the author's knowledge, a new concept and based on the reformulation of the cubic Klein-Gordon equation as a system of coupled nonlinear Helmholtz equations involving suitable conditions on the far field behavior. 
\end{abstract}

\maketitle
 
\allowdisplaybreaks

\section{Introduction and Main Results}

\smallskip

We construct real-valued solutions $U(t, x)$ of the cubic Klein-Gordon equation
\begin{equation}\label{eq_wave}
	\partial_t^2 U - \Delta U + m^2 U = \Gamma(x) \: U^3
	\qquad
	\text{on } \R \times \R^3
\end{equation}
where $\Gamma \in L^\infty_\mathrm{rad}(\R^3) \cap C^1_\mathrm{loc}(\R^3)$ and $m > 0$ is a (mass) parameter. Here we restrict ourselves to the case of three space dimensions which is the most relevant one for applications in physics and which allows to use the tools established in~\cite{own_cubic}. Throughout, the notations $\partial_{1, 2, 3}, \nabla, \Delta, D^2$ refer to differential operators acting on the space variables.
The solutions we aim to construct are polychromatic, that is, they take the form
\begin{align}\label{eq_poly}
	&U(t, x) = u_0(x) + \sum_{k=1}^\infty  2 \,\cos(\omega k t)  u_k(x) = \sum_{k \in \Z}  \mathrm{e}^{\i \omega k t}  u_k(x)
	\\
	\nonumber
	&\text{where} \: \:
	u_k \in X_1= \left\{ u \in C_\mathrm{rad}(\R^3, \R) \big| \: \norm{(1 + |\cdot|^2)^\frac{1}{2} u}_\infty < \infty \right\} 
	\subseteq L^4_\mathrm{rad}(\R^3),
	\: \: u_{-k} = u_k
	\\
	\nonumber
	&\text{and (for simplicity)} \: \:  \omega > m.
\end{align}
Such solutions are periodic in time and localized as well as radially symmetric in space. They are sometimes referred to as breather solutions, c.f. the ``Sine-Gordon breather'' in~\cite{Ablowitz}, equation (28). 
The construction of breather solutions is of particular interest since, as indicated in a study~\cite{Weinstein} on perturbations of the Sine-Gordon breather, Birnir, McKean and Weinstein conjecture that  ``for the general nonlinear wave equation [author's note: in 1+1 dimensions], breathing [...] takes place only for isolated nonlinearities'', see ~\cite[p.1044]{Weinstein}. This conjecture is supported by recent existence results for breathers for the 1+1 dimensional wave equation with specific, carefully designed potentials which we comment on below. Our results, however, indicate that the situation might be entirely different for weakly localized breathers for the Klein-Gordon equation in 1+3 dimensions, in the sense that such breather solutions are abundant even in ``simple'' settings.

\smallskip

We will find breather solutions of~\eqref{eq_wave} with $u_k \not\equiv 0$ for at least two distinct integers $k \in \N_0$ by rewriting it into an infinite system of stationary equations for the functions $u_k$. Indeed, inserting~\eqref{eq_poly}, a short and formal calculation leads to 
\begin{subequations}\label{eq_stationary_0}
\begin{align}
	- \Delta u_0 + m^2 \, u_0 &= \Gamma(x) \: \left( \b{u} \star \b{u} \star \b{u} \right)_0,
	\label{eq_stationary_1}
	\\
	- \Delta u_{k} - (\omega^2 k^2 - m^2) u_{k} &= \Gamma(x) \: \left( \b{u} \star \b{u} \star \b{u} \right)_{k}
	\label{eq_stationary_2}
	\qquad \text{for } k \in \Z \setminus \{ 0 \}.
\end{align}
\end{subequations}
In fact, \eqref{eq_stationary_2} includes \eqref{eq_stationary_1}, but we intend to separate the ``Schrödinger'' equation characterized by $0 \not\in \sigma(- \Delta + m^2)$ from the infinite number of ``Helmholtz'' equations characterized by $0 \in \sigma(- \Delta - (\omega^2 k^2 - m^2))$, $k \neq 0$. Our construction of breathers for~\eqref{eq_wave} relies on new methods for such Helmholtz equations introduced in~\cite{own_cubic} which exploit the so-called far field properties of their solutions and lead to a rich bifurcation structure. These methods will be sketched only briefly in the main body of this paper; more details will be given in Section~\ref{sect_helmholtz} at the end (which can be read independently).

\smallskip

The solutions we obtain bifurcate from any given stationary (radial) solution
 $w_0 \in X_1$, $w_0 \not \equiv 0$ of the Klein-Gordon equation~\eqref{eq_wave}. That is, $w_0$ solves the stationary nonlinear Schrödinger equation
\begin{align}\label{eq_w0}
	- \Delta w_0 + m^2 \, w_0 = \Gamma (x) \: w_0^3
	\quad \text{on }\R^3;
\end{align}
regarding existence of such $w_0$, cf. Remark~\ref{rmk}~(b). Let us remark briefly that all (distributional) solutions of~\eqref{eq_w0} in $X_1 \subseteq L^4_\text{rad}(\R^3)$ are twice differentiable by elliptic regularity. 
In order to make bifurcation theory work, we impose the following nondegeneracy  assumption: 
\begin{align}\label{eq_nondegenerate,S}
	q_0 \in X_1, \: \: - \Delta q_0 + q_0 = 3 \Gamma(x) \: w_0^2 \: q_0
	\text{ on } \R^3
	\qquad \text{implies } 
	\qquad
	q_0 \equiv 0.
\end{align}
We comment on this assumption in Remark~\ref{rmk}~(c) below. 
In particular,~\eqref{eq_nondegenerate,S} and our main result presented next hold if $\Gamma$ is constant and $w_0$ is a (positive) ground state of~\eqref{eq_w0}. We now present our main result.

\begin{thm}\label{thm_poly}
Let $\Gamma \in L^\infty_\mathrm{rad}(\R^3) \cap C^1_\mathrm{loc}(\R^3)$, $\omega > m > 0$ and assume there is some stationary solution $U^0(t,x) = w_0(x)$, $w_0 \not \equiv 0$ of the cubic Klein-Gordon equation~\eqref{eq_wave}, i.e. $w_0 \in X_1$ solving~\eqref{eq_w0}. Assume further that $w_0$ is nondegenerate in the sense of~\eqref{eq_nondegenerate,S}. Then for every $s \in \N$  there exist an open interval $J_s \subseteq \R$ with $0 \in J_s$ and a family $(U^{\alpha})_{\alpha \in J_s} \subseteq C^2(\R, X_1)$ with the following properties:
\begin{itemize}
\item[(i)] All $U^\alpha$ are time-periodic, twice continuously differentiable classical solutions of~\eqref{eq_wave} of the polychromatic form~\eqref{eq_poly},
\begin{align*}
	U^\alpha(t, x) = u_0^\alpha(x) + \sum_{k=1}^\infty  2 \,\cos(\omega k t)  u_{k}^\alpha(x).
\end{align*}
\item[(ii)] The map $\alpha \mapsto (u_k^\alpha)_{k \in \N_0}$ is smooth in the topology of $\ell^1(\N_0, X_1)$ with
\begin{align*}
	\frac{\mathrm{d}}{\mathrm{d}\alpha}\bigg|_{\alpha = 0} u_k^\alpha \not\equiv 0
	\quad \text{if and only if} \quad k = s
\end{align*}
(``excitation of the s-th mode'').
In particular, for sufficiently small $\alpha \neq 0$, these solutions are non-stationary.
Moreover, for different values of $s$,  the families of solutions mutually differ close to $U^0$.
\item[(iii)] If we assume additionally $\Gamma(x) \neq 0$ for almost all $x \in \R^3$, then every nonstationary polychromatic solution $U^\alpha$ possesses infinitely many nonvanishing modes $u_k^\alpha$.
\end{itemize}
\end{thm}

\begin{rmk}\label{rmk}
\begin{itemize}
\item[(a)]
We require continuity of $\Gamma$ since we use the functional analytic framework of~\cite{own_cubic}. The existence and continuity of $\nabla \Gamma$ will be exploited in proving that $U^\alpha$ is twice differentiable. This assumption as well as $\Gamma \neq 0$ almost everywhere in (iii) might be relaxed; however, this study does not aim at the most general setting for the coefficients but rather focuses on the introduction of the setup for the existence result.
\item[(b)]
The existence of stationary solutions of the Klein-Gordon equation~\eqref{eq_wave} resp. of solutions to~\eqref{eq_w0} can be guaranteed under additional assumptions on $\Gamma$. We refer to \cite{Lions}, Theorem~I.2 and Remarks~I.5,~I.6 by Lions for positive (ground state) solutions and to Theorems~2.1 of~\cite{Willem1},~\cite{Willem2} by Bartsch and Willem for bound states.
\item[(c)]
In some special cases, nondegeneracy properties like~\eqref{eq_nondegenerate,S} have been verified, e.g. by Bates and Shi~\cite{bates} in Theorem~5.4~(6), or by Wei~\cite{wei} in Lemma~4.1, both assuming that $w_0$ is a ground state solution of~\eqref{eq_w0} in the autonomous case with constant positive $\Gamma$. It should be pointed out that, although the quoted results discuss nondegeneracy in a setting on the Hilbert space $H^1(\R^3)$, the statements can be adapted to the topology of $X_1$, as we will demonstrate  in Lemma~\ref{lem_nondegeneracy}.
\item[(d)]
The assumption $\omega > m$ on the frequency ensures that the stationary system~\eqref{eq_stationary_0} contains only one equation of Schrödinger type. This avoids further nondegeneracy assumptions on higher modes, which would not be covered by the previously mentioned results in the literature.
\item[(e)]
The above result provides, locally, a multitude of families of breathers bifurcating from every given stationary solution characterized by different values of $s$, $\omega$ and possibly certain asymptotic parameters, see Remark~\ref{rmk_previouschapter} below. \\ It would be natural, further, to ask for the global bifurcation picture given some trivial family $\mathcal{T} = \{ (w_0, \lambda) \, | \, \lambda \in \R \}$. (Here $\lambda \in \R$ denotes a bifurcation parameter which in our case is not visible in the differential equation and thus will be properly introduced later.)  Typically, global bifurcation theorems state that a maximal bifurcating continuum of solutions $(U, \lambda)$ emanating from $\mathcal{T}$ at $(w_0, \lambda_0)$ is unbounded unless it returns to $\mathcal{T}$ at some point $(w_0, \lambda_0')$, $\lambda_0' \neq \lambda_0$. In the former (desirable) case, however, a satisfactory characterization of global bifurcation structures should provide a criterion whether or not unboundedness results from another stationary solution $w_1 \neq w_0$ with $\{ (w_1, \lambda) \, | \, \lambda \in \R \}$ belonging to the maximal continuum. Since it is not obvious at all whether and how such a criterion might be derived within our framework, we focus on the local result, which already adds new aspects to the state of knowledge about the existence of breather solutions summarized next.
\end{itemize}
\end{rmk}

\subsection{An Overview of Literature}

\subsubsection*{Polychromatic Solutions}

The results in Theorem~\ref{thm_poly} can and should be compared with recent findings on breather (that is to say, time-periodic and spatially localized) solutions of the wave equation with periodic potentials $V(x), q(x) = c \cdot V(x) \geq 0$,
\begin{align}\label{eq_wave1D}
	V(x) \partial_t^2 U - \partial_x^2 U + q(x) U = \Gamma(x) U^3
	\qquad  \text{on } \R \times \R.
\end{align}
Such breather solutions have been constructed by Schneider et al., see Theorem~1.1~in~\cite{Schneider}, and Hirsch and Reichel, see Theorem~1.3~in~\cite{Hirsch}, respectively. 
In brief, the main differences to the results in this article are that the authors of~\cite{Schneider},~\cite{Hirsch} consider a setting in one space dimension and obtain strongly spatially localized solutions, which requires a comparably huge technical effort. We give some details:
Both existence results are established using a polychromatic ansatz, which reduces the time-dependent equation to an infinite set of stationary problems with periodic coefficients, see \cite{Schneider},~p.~823, resp. \cite{Hirsch},~equation~(1.2).  
The authors of~\cite{Schneider} apply spatial dynamics and center manifold reduction; their ansatz is based on a very explicit choice of the coefficients $q, V, \Gamma$. The approach in~\cite{Hirsch} incorporates more general potentials and nonlinearities and is based on variational techniques. It provides ground state solutions, which are possibly ``large'' - in contrast to our local bifurcation methods, which only yield solutions close to a given stationary one as described in Theorem~\ref{thm_poly}, i.e. with a typically ``small'' time-dependent contribution. 
\\
Periodicity of the potentials in~\eqref{eq_wave1D} is explicitly required since it leads to the occurrence of spectral gaps when analyzing the associated differential operators of the stationary equations. 
In contrast to the Helmholtz methods introduced here, the authors both of~\cite{Schneider} and of~\cite{Hirsch} strive to construct the potentials in such way that  $0$ lies in the aforementioned spectral gaps, and moreover that the distance between $0$ and the spectra has a positive lower bound. This is realized by assuming a certain ``roughness'' of the potentials, referring to the step potential defined in Theorem~1.1~of~\cite{Schneider} and to the assumptions (P1)-(P3) in~\cite{Hirsch} which allow potentials with periodic spikes modeled by Dirac delta distributions, periodic step potentials or some specific, non-explicit potentials in $H^r_\text{rad}(\R)$ with $1 \leq r < \frac{3}{2}$ (see~\cite{Hirsch}, Lemma~2.8). 

\smallskip

Let us summarize that the methods for constructing breather solutions of~\eqref{eq_wave1D} outlined above can handle periodic potentials but require irregularity, are very restrictive concerning the form of the potentials and involve a huge technical effort in analyzing spectral properties based on Floquet-Bloch theory. The Helmholtz ansatz presented in this article provides a technically elegant and short approach suitable for constant potentials; in the context of breather solutions, it is new in the sense that it provides breathers with slow decay, it provides breathers on the full space $\R^3$, and it provides breathers for simple (constant) potentials.

\subsubsection*{The Klein-Gordon Equation as a Cauchy Problem}

Possibly due to its relevance in physics, there is a number of classical results in the literature concerning the nonlinear Klein-Gordon equation. 
The fundamental difference to the results in this article is that the vast majority of these concerns the Cauchy problem of the Klein-Gordon equation, i.e.
\begin{equation}\label{eq_NKG}
\begin{split}
	&\partial_t^2 U - \Delta U + m^2 \,  U = \pm U^3
	\quad \text{on } [0, \infty) \times \R^3
	\\
	&U(0, x) = f(x), \: \partial_t U(0, x) = g(x)
	\quad \text{ on } \R^3
\end{split}
\end{equation}
for suitable initial data $f, g: \R^3 \to \R$. Usually, the dependence of the nonlinearity on $U$ is much more general (allowing also derivatives of $U$) and the space dimension is not restricted to $N=3$. On the other hand, most results in the literature only concern the autonomous case, which is why we set in this discussion $\Gamma \equiv \pm 1$. 

An overview of the state of knowledge towards the end of the 1970s can be found e.g. in~\cite{StraussNKG} by Strauss, who discusses among other topics global existence (Theorem~1.1), regularity and uniqueness (Theorem~1.2), blow-up (Theorem~1.4) and scattering
%, i.e. convergence to solutions of the free Klein-Gordon equation as $t \to \infty$ 
(Theorem~4.1). 
In the first-mentioned result, which is originally due to Jörgens, global existence of distributional solutions with locally as well as globally finite energy 
\begin{align*}
	E_B[U(t,\,\cdot\,)]
	= \frac{1}{2} \int_B |\partial_t U(t,x)|^2 + |\nabla U(t,x)|^2 + m^2 |U(t,x)|^2 \: \mathrm{d}x
	+ \frac{1}{4} \int_B |U(t, x)|^4 \: \mathrm{d}x, 
	\:\:
	B \subseteq \R^3
\end{align*}
is proved provided $\Gamma \equiv -1$. Following a classical strategy for evolution problems, local existence is shown by means of a fixed point iteration, and global existence can be obtained by an iteration argument based on energy conservation. For $\Gamma \equiv +1$, Theorem~1.4 due to Keller and Levine demonstrates the existence of blow-up solutions.
\\
During the following decade, Klainerman~\cite{Klainerman1, Klainerman2} and Shatah~\cite{Shatah1, Shatah2} independently developed new techniques leading to significant improvements in the study of uniqueness questions and of the asymptotic behavior of solutions as $t \to \infty$. These results work in settings with high regularity and admit more general nonlinearities with growth assumptions for small arguments, which includes the cubic one as a special case. In particular, Klainerman and Shatah prove the convergence to solutions of the free Klein-Gordon equation and show uniform decay rates of solutions as $t \to \infty$. In the case of a cubic nonlinearity, these results only apply if the space dimension is at least $2$. This is why, more recently, the question of corresponding uniqueness and convergence properties for cubic nonlinearities in $N = 1$ space dimensions has attracted attention; we wish to mention at least some of the related papers. For explicit choices of the cubic nonlinearity, there are results by Moriyama and by Delort, see Theorem~1.1~of~\cite{Moriyama} resp. Th\'{e}or\`{e}mes~1.2,~1.3~in~\cite{Delort}. Only the latter result allows a nonlinearity of the form $\pm U^3$ not containing derivatives (see~\cite{Delort}, Remarque~1.4); however, the initial data are assumed to have compact support. 
Global existence, uniqueness, decay rates and scattering exclusively for the nonlinearity $\pm U^3$ can be found in Corollary~1.2~of~\cite{Hayashi} by Hayashi and Naumkin.
 
The relation to our results is not straightforward since the bifurcation methods automatically provide solutions $U^\alpha$ which exist globally in time irrespective of the sign (or even of a possible $x$-dependence) of $\Gamma$ and which do not decay as $t \to \infty$, and there is no special emphasis on the role of the initial values $U^\alpha(0, x), \nabla U^\alpha(0, x)$ along the bifurcating branches. 
Our methods instead focus on several global properties of the solutions $U^\alpha(t, x)$ such as periodicity in time and localization as well as decay rates in space, i.e. the defining properties of breathers.

\subsection{Research Perspectives}

Apart from bifurcation methods, nonlinear Helmholtz equations and systems can also be discussed in a ``dual'' variational framework as introduced by Ev\'{e}quoz and Weth~\cite{EvequozWeth}. This might offer another way to analyze the system~\eqref{eq_stationary_0} leading  to ``large'' breathers in the sense that they are not close to a given stationary solution as the ones constructed in Theorem~\ref{thm_poly}. Furthermore, such an ansatz might be a promising step towards extensions to non-constant, e.g. periodic potentials.

\section{The Proof of Theorem~\ref{thm_poly}}

\subsection{The Functional-Analytic Setting}

We look for polychromatic solutions as in~\eqref{eq_poly} with coefficients $\b{u} = (u_k)_{k \in \Z} \in \mathcal{X}_1$ where
\begin{align*}
	\mathcal{X}_1 := \ell^1_\text{sym}(\Z, X_1)
	:= \left\{
	(u_k)_{k \in \Z} \: \bigg| \: u_k = u_{-k} \in X_1, \norm{(u_k)_{k \in \Z}}_{\mathcal{X}_1} := \sum_{k \in \Z} \norm{u_k}_{X_1} < \infty
	\right\}.
\end{align*}
The Banach space $X_1$ has been defined in~\eqref{eq_poly}; it prescribes a decay rate which is the natural one for solutions of Helmholtz equations as in~\eqref{eq_stationary_2}, see also Section~\ref{sect_helmholtz}.
Throughout, we denote by $\b{w} = (\delta_{k, 0} w_0)_{k \in \Z} = ( ..., 0, w_0, 0, ...) $ the stationary solution with $w_0 \in X_1 \cap C^2_\mathrm{loc}(\R^3)$ fixed according to equation~\eqref{eq_w0}. We will find polychromatic solutions of~\eqref{eq_wave} by solving the countably infinite Schrödinger-Helmholtz system~\eqref{eq_stationary_1},~\eqref{eq_stationary_2}, which is equivalent to~\eqref{eq_wave},~\eqref{eq_poly} on a formal level; for details including convergence of the polychromatic sum in~\eqref{eq_poly}, see Proposition~\ref{prop_mapping}. 

\medskip

Our strategy is then as follows:
\begin{itemize}
\item[$\triangleright$]
Intending to apply bifurcation techniques, we have to analyze the linearized version of the infinite-dimensional system~\eqref{eq_stationary_1},~\eqref{eq_stationary_2}, which resembles the one of the two-component system discussed by the author in~\cite{own_cubic}. 
We therefore summarize, for the reader's convenience, a collection of results concerning the linearized setting in Proposition~\ref{prop_chapter2}. 
\item[$\triangleright$]
We then present a suitable setup for bifurcation theory; in particular, we introduce a bifurcation parameter which is not visible in the differential equation but appears in the so-called far field of the functions $u_k$, more specifically a phase parameter in the leading-order contribution as $|x| \to \infty$. 
\item[$\triangleright$]
The aforementioned fact that solutions of~\eqref{eq_stationary_1},~\eqref{eq_stationary_2} obtained in this setting provide polychromatic, classical solutions of the Klein-Gordon equation~\eqref{eq_wave} will be proved as a part of Proposition~\ref{prop_mapping} below. Indeed, regarding differentiability, we will see that the choice of suitable asymptotic conditions will ensure uniform convergence and hence smoothness properties of the infinite sums defining the polychromatic states.
\item[$\triangleright$]
Finally, in Proposition~\ref{prop_kernel}, we essentially verify the assumptions of the Crandall-Rabinowitz Bifurcation Theorem. 
\end{itemize}
After that, we are able to give a very short proof of Theorem~\ref{thm_poly}.
The auxiliary results will be proved in Section~\ref{ch_wave-proofs}. The final Section~\ref{sect_helmholtz} provides some more details on the theory of linear Helmholtz equations in $X_1$.

\medskip

Throughout, we denote the convolution in $\R^3$ by the symbol $\ast$ and use $\star$ in the convolution algebra $\ell^1$. Extending the notation defined above, for $q \geq 0$, we let
\begin{align*}
	&X_q := \left\{ u \in C_\text{rad}(\R^3, \R) \, | \, \norm{u}_{X_q} < \infty \right\}
	&& \text{with } \norm{u}_{X_q} := \sup_{x \in \R^3} (1 + |x|^2)^{q/2} |u(x)|,
	\\
	& \mathcal{X}_q := \ell^1_\text{sym}(\Z, X_q)
	&& \text{with } \norm{\b{u}}_{\mathcal{X}_q} := \norm{(u_k)_k}_{\mathcal{X}_q} := \sum_{k \in \Z} \norm{u_k}_{X_q}.
\end{align*}

\begin{prop}\label{prop_convolution-x1}
The convolution of sequences $\b{u}^{(1)}, \b{u}^{(2)}, \b{u}^{(3)} \in \mathcal{X}_1$ is well-defined in a pointwise sense and satisfies $\b{u}^{(1)} \star \b{u}^{(2)} \star \b{u}^{(3)} \in \mathcal{X}_3$. Moreover, we have the estimate
\begin{align*}
	\norm{\b{u}^{(1)} \star \b{u}^{(2)} \star \b{u}^{(3)}}_{\mathcal{X}_3}
	\leq \norm{\b{u}^{(1)}}_{\mathcal{X}_1} \norm{\b{u}^{(2)}}_{\mathcal{X}_1} \norm{\b{u}^{(3)}}_{\mathcal{X}_1}.
\end{align*} 
\end{prop}

We rewrite the system~\eqref{eq_stationary_1},~\eqref{eq_stationary_2} using $\b{u} = \b{w} + \b{v}$ with $\b{w} = ( ..., 0, w_0, 0, ...) $; then,
\begin{align}\label{eq_stationary-v}
	- \Delta v_k - (\omega^2 k^2 - m^2) \, v_k = \Gamma (x) \cdot  
	\left[ \left((\b{w} + \b{v}) \star (\b{w} + \b{v}) \star (\b{w} + \b{v})\right)_k - \delta_{k, 0} w_0^3 \right]
	\quad \text{on }\R^3.
\end{align}
We will find solutions of this system of differential equations by solving instead a system of coupled convolution equations which, for $k \not\in \{ 0,  \pm s \}$, have the form $v_k = \mathcal{R}_{\mu_k}^{\tau_k} [f_k]$. Here $f_k$ represents the right-hand side of~\eqref{eq_stationary-v}, $\mu_k := \omega^2 k^2 - m^2$, and the coefficients $\tau_k \in (0, \pi)$ will have to be chosen properly according to a nondegeneracy condition. The convolution operators 
$$\mathcal{R}_{\mu}^{\tau} = \frac{\sin(|\,\cdot\,| \sqrt{\mu} + \tau)}{4\pi \sin(\tau) |\,\cdot\,|} \: \ast \: : X_3 \to X_1
\qquad (\mu > 0, \: \: 0 < \tau < \pi)$$
can be viewed as resolvent-type operators for the Helmholtz equation $(- \Delta - \mu) v = f$ on $\R^3$ involving an asymptotic condition on the far field of the solution $v$, namely 
$$ |x| \: v(x) \sim \sin(|x| \sqrt{\mu} + \tau) + O\left(\frac{1}{|x|}\right) 
\quad \text{as} \quad |x| \to \infty. $$ 
Such conditions are required since the homogeneous Helmholtz equation $(- \Delta - \mu) v = 0$ has smooth nontrivial solutions in $X_1$ (known as Herglotz waves), which are all multiples of 
$$ \tilde{\Psi}_\mu (x) := \frac{\sin(|x| \sqrt{\mu})}{4 \pi |x|} \quad (x \neq 0). $$
We refer to Section~\ref{sect_helmholtz}, more precisely Lemma~\ref{lem_linH}, for details; the case $\tau = 0$ requires a larger technical effort and is presented in Lemma~\ref{lem_linH-0}. This involves linear functionals $\alpha^{(\mu)}, \beta^{(\mu)} \in X_1'$ which, essentially, yield the coefficients of the sine resp. cosine terms in the asymptotic expansion above.  Relying on these tools and notations, we summarize the relevant facts on the linearized versions of the Helmholtz equations~\eqref{eq_stationary_2} in the following Proposition. 

\begin{prop}\label{prop_chapter2}
Let $w_0 \in X_1$ be a solution of equation~\eqref{eq_w0} with $\Gamma \in L^\infty_\mathrm{rad}(\R^3) \cap C_\mathrm{loc}(\R^3)$ and $\omega > m > 0$; define $\mu_k := \omega^2 k^2 - m^2$. For every $k \in \Z \setminus \{ 0 \}$, there exists (up to a multiplicative constant) a unique nontrivial and radially symmetric solution $q_k \in X_1$ of 
\begin{subequations}\label{eq_qk}
\begin{align}\label{eq_qk-1}
	- \Delta q_k - \mu_k \, q_k = 3 \, \Gamma(x) w_0^2(x) \: q_k \qquad \text{on }\R^3.
\end{align} 
It is twice continuously differentiable and satisfies, for some $c_k \neq 0$ and $\sigma_k \in [0, \pi)$,
\begin{align}\label{eq_qk-2}
	q_k(x) = c_k \cdot \frac{\sin(|x| \, \sqrt{\mu_k} + \sigma_k)}{|x|} + O\left(\frac{1}{|x|^2}\right)
	\quad \text{as } |x| \to \infty.
\end{align}
\end{subequations}
The equations~\eqref{eq_qk-1},~\eqref{eq_qk-2} are equivalent to the convolution identities
\begin{align*}
	\begin{cases}
	q_k	= 
	3 \:  \mathcal{R}_{\mu_k}^{\sigma_{k}} [\Gamma w_0^2 \, q_{k}]
	= 
	3 \: \left( \mathcal{R}_{\mu_k} [\Gamma w_0^2 \, q_{k}] 
	 + \cot(\sigma_k) \tilde{\mathcal{R}}_{\mu_k} [\Gamma w_0^2 \, q_{k}] \right)
	& \text{if } \sigma_k \in (0, \pi),
	\\
	q_k = 3 \:  \mathcal{R}_{\mu_k}^{\pi/2} [\Gamma w_0^2 \, q_{k}] 
	+ \left( \alpha^{(\mu_k)}(q_{k}) + \beta^{(\mu_k)}(q_{k}) \right) \cdot \tilde{\Psi}_{\mu_k}
	& \text{if } \sigma_k = 0.
	\end{cases}
\end{align*}
For all $k \in \Z$, $\cos(\sigma_k) \, \beta^{(\mu_k)}(q_{k}) = \sin(\sigma_k) \, \alpha^{(\mu_k)}(q_{k})$.
\end{prop}

The existence statement and the asymptotic properties in~\eqref{eq_qk} can be proved using the Prüfer transformation, see~\cite{own_cubic}, Proposition~6; the statements in the second part are consequences of Lemmas~\ref{lem_linH}~and~\ref{lem_linH-0} in the final Section~\ref{sect_helmholtz}. For these results to apply we have assumed initially that $\Gamma$ is continuous and bounded, whence $3 \, \Gamma w_0^2 \in X_2$.

We now present the general assumptions valid throughout the following construction and the proof of Theorem~\ref{thm_poly}. We let $\sigma_k$ for $k \in \Z \setminus \{ 0 \}$ as in Proposition~\ref{prop_chapter2} above and fix $s \in \N$, recalling that we aim to ``excite the $s$-th mode'' in the sense of Theorem~\ref{thm_poly}~(ii). With this, let us introduce 
\begin{equation}\label{eq_assumptions}
	\tau_{\pm s} := \sigma_{\pm s},
	\qquad
	\tau_k := \begin{cases}
		\frac{\pi}{4}		& \text{if } \sigma_k \neq \frac{\pi}{4}, 
		\\
		\frac{3\pi}{4}		& \text{if } \sigma_k = \frac{\pi}{4}
	\end{cases}
	\quad \text{  for } k \in \Z \setminus \{ 0, \pm s \},
\end{equation}
see also Remark~\ref{rmk_previouschapter}~(b).
Thus in particular $\tau_k \neq \sigma_k$ for $k \in \Z \setminus \{ 0, \pm s \}$, and we conclude from the uniqueness statement in Proposition~\ref{prop_chapter2} the nondegeneracy property 
\begin{subequations}\label{eq_nondegenerate}
\begin{equation}\label{eq_nondegenerate,1}
	k \in \Z \setminus \{0, \pm s \}, \quad q \in X_1, \quad 
	q	= 3 \: \mathcal{R}_{\mu_k}^{\tau_k} [\Gamma w_0^2 \, q]
	\qquad \Rightarrow \qquad
	q \equiv 0;
\end{equation}
for the $0$-th mode, using the resolvent $\mathcal{P}_{\mu_0} = (- \Delta + \mu_0)^{-1}: X_3 \to X_1$ (see Lemma~\ref{lem_linS}), the corresponding property is assumed in~\eqref{eq_nondegenerate,S}:
\begin{equation}\label{eq_nondegenerate,2}
	q \in X_1, \quad 
	q	= 3 \: \mathcal{P}_{\mu_0} [\Gamma w_0^2 \, q]
	\qquad \Rightarrow \qquad
	q \equiv 0.
\end{equation}
\end{subequations}
We now introduce a map the zeros of which provide solutions of the system~\eqref{eq_stationary-v}. 
Throughout, we use the shorthand notation $\b{u} = \b{v} + \b{w}$ for $\b{v} \in \mathcal{X}_1$ and the stationary solution $\b{w} = (..., 0, w_0, 0, ...)$.
As above, we have to distinguish the cases $\tau_s \in (0, \pi)$ and $\tau_s = 0$. (In the following, please recall that we consider some fixed $s \neq 0$.)
For $0 < \tau_{\pm s} < \pi$, we introduce $F: \: \mathcal{X}_1 \times \R \to \mathcal{X}_1$ via
\begin{subequations}\label{eq_FG}
\begin{equation}\label{eq_F}
F(\b{v}, \lambda)_{k}
:= v_k -
\begin{cases}
	\mathcal{P}_{\mu_0}
	\left[ \Gamma \: \left(\b{u} \star \b{u} \star \b{u}\right)_0 
	- \Gamma \: w_0^3  \right] 
	& k = 0,
	\\
	\mathcal{R}_{\mu_s}^{\pi/2}
	\left[ \Gamma \: \left(\b{u} \star \b{u} \star \b{u}\right)_{\pm s} \right]
	\\
	\qquad 
	+ (\cot(\tau_{\pm s}) - \lambda) \tilde{\Psi}_{\mu_s} \ast
	\left[ \Gamma \: \left(\b{u} \star \b{u} \star \b{u}\right)_{\pm s} \right]
	& k = \pm s,
	\\
	\mathcal{R}_{\mu_k}^{\tau_{k}}
	\left[ \Gamma \: \left(\b{u} \star \b{u} \star \b{u}\right)_{k} \right]
	& \text{else}.
\end{cases}
\end{equation}
Similarly, if $\sigma_s = 0$, we define $G: \: \mathcal{X}_1 \times \R \to \mathcal{X}_1$ by
\begin{equation}\label{eq_G}
G(\b{v}, \lambda)_{k}
:= v_k - 
\begin{cases}
	\mathcal{P}_{\mu_0}
	\left[ \Gamma \: \left(\b{u} \star \b{u} \star \b{u}\right)_0 
	-  \Gamma \: w_0^3  \right] 
	& k = 0,
	\\
	\mathcal{R}_{\mu_s}^{\pi/2}
	\left[ \Gamma \: \left(\b{u} \star \b{u} \star \b{u}\right)_{\pm s} \right] 
	\\
	\qquad
	+ (1 - \lambda)
	\left( \alpha^{(\mu_s)}(v_{\pm s}) + \beta^{(\mu_s)}(v_{\pm s}) \right) \tilde{\Psi}_{\mu_s}
	& k = \pm s, 
	\\
	\mathcal{R}_{\mu_k}^{\tau_{k}}
	\left[ \Gamma \: \left(\b{u} \star \b{u} \star \b{u}\right)_{k} \right]
	& \text{else}.
\end{cases}
\end{equation}
\end{subequations}

The following result collects some basic properties of the maps $F$ and $G$ and the polychromatic states related to their zeros. 

\begin{prop}\label{prop_mapping}
Let $s \in \N$ and $(\tau_k)_{k \in \Z}$ be chosen as in~\eqref{eq_assumptions}. 
The maps $F, G: \mathcal{X}_1 \times \R \to \mathcal{X}_1$ are well-defined and smooth with $F(\b{0}, \lambda) = G(\b{0}, \lambda) = \b{0}$ for all $\lambda \in \R$.
Further, if $F(\b{v}, \lambda) = \b{0}$ resp. $G(\b{v}, \lambda) = \b{0}$ for some $ \b{v} \in \mathcal{X}_1, \lambda \in \R$, then $\b{v}$ solves the stationary system~\eqref{eq_stationary-v} and
\begin{align*}
	U(t, x) := w_0(x) + v_0(x)  + \sum_{k=1}^\infty  2 \,\cos(\omega k t)  v_k(x)
	\qquad
	(t \in \R, x \in \R^3)
\end{align*}
defines a twice continuously differentiable, classical solution $U \in C^2(\R, X_1)$ of the Klein-Gordon equation~\eqref{eq_wave}.
\end{prop}
Again, the proof can be found in Section~\ref{ch_wave-proofs}. We will even show that $U \in C^\infty(\R, X_1)$. For the derivatives of $F$ resp. $G$ with respect to the Banach space component $\b{v} \in \mathcal{X}_1$, we will verify  the following explicit formulas: Letting $\b{q} \in \mathcal{X}_1$ and abbreviating $\b{u} := \b{v} + \b{w}$, 
\begin{subequations}\label{eq_DFG}
\begin{align}\label{eq_DF}
	(DF(\b{v}, \lambda)[\b{q}])_k
	= q_k - \begin{cases}
	3 \: \mathcal{P}_{\mu_0} [\Gamma  \left(\b{q} \star \b{u} \star \b{u}\right)_0]
	& k = 0,
		\\
	3 \: \mathcal{R}_{\mu_s}^{\pi/2}
	\left[ \Gamma  \left(\b{q} \star \b{u} \star \b{u}\right)_{\pm s} \right]
	& 
	\\
	\:
	+ 3 \: (\cot(\tau_s) - \lambda) \tilde{\Psi}_{\mu_s} \ast
	\left[ \Gamma  \left(\b{q} \star \b{u} \star \b{u} \right)_{\pm s} \right]
	& k = \pm s,
	\\
	3 \: \mathcal{R}_{\mu_k}^{\tau_k}
	\left[ \Gamma  \left(\b{q} \star \b{u} \star \b{u}\right)_k \right]
	& \text{else};
	\end{cases}
\end{align}
\begin{align}\label{eq_DG}
	(DG(\b{v}, \lambda)[\b{q}])_k
	= q_k - \begin{cases}
	3 \: \mathcal{P}_{\mu_0} [\Gamma  \left(\b{q} \star \b{u} \star \b{u}\right)_0]
	& k = 0,
	\\
	3 \: \mathcal{R}_{\mu_s}^{\pi/2}
	\left[ \Gamma  \left(\b{q} \star \b{u} \star \b{u}\right)_{\pm s} \right]
	& 
	\\
	\:
	+  (1 - \lambda) \!
	\left( \alpha^{(\mu_s)}(q_{\pm s}) + \beta^{(\mu_s)}(q_{\pm s}) \right) \! \tilde{\Psi}_{\mu_s}
	& k = \pm s,
		\\
	3 \: \mathcal{R}_{\mu_k}^{\tau_k}
	\left[ \Gamma  \left(\b{q} \star \b{u} \star \b{u}\right)_k \right]
	& \text{else}.
	\end{cases}
\end{align}
\end{subequations}

\begin{rmk}\label{rmk_previouschapter}
\begin{itemize}
\item[(a)]
As earlier announced, we now see that the bifurcation parameter $\lambda$ appears only in the asymptotic expansions of the $s$-th components $v_{\pm s}$ of the solutions and not in the differential equation~\eqref{eq_wave}. This is different from~\cite{own_cubic} where the bifurcation parameter takes the role of a coupling parameter of the Helmholtz system. 
\item[(b)]
The choice of the parameters $\tau_k$ in equation~\eqref{eq_assumptions} is far from unique. Indeed, one could instead consider any configuration satisfying
\begin{align*}
	\tau_k = \tau_{-k} \neq \sigma_k \: \text{ for all } k \in \Z \setminus \{ \pm s \},
	\qquad
	\{ \tau_k \, | \, k \in \Z \setminus \{ \pm s \} \} \subseteq (\delta, \pi - \delta)
\end{align*}
for some $\delta \in \left(0, \frac{\pi}{2} \right)$. The former condition is required for the nondegeneracy statement~\eqref{eq_nondegenerate,1}, and the latter will be used to obtain uniform decay estimates in the proof of Proposition~\ref{prop_mapping}, see Lemma~\ref{lem_uniformdecay}.
\\
However, as in~\cite{own_cubic}, the question whether another choice of $\tau_k$ leads to different bifurcating families is still open.
Hence we discuss only the explicit choice in~\eqref{eq_assumptions}.
\end{itemize}
\end{rmk}

In the so-established framework, we intend to apply the Crandall-Rabinowitz Bifurcation Theorem. The next result shows that its assumptions are satisfied.

\begin{prop}[Simplicity and transversality]\label{prop_kernel}
Let $s \in \N$ and $(\tau_k)_{k \in \Z}$ be chosen as in~\eqref{eq_assumptions}. 
The linear operator $DF(\b{0}, 0): \mathcal{X}_1 \to \mathcal{X}_1$ is 1-1-Fredholm with a kernel of the form
\begin{align*}
	\ker DF(\b{0}, 0) = \mathrm{span } \, \{ \b{q} \}
	\qquad \text{where } q_k \neq 0 \text{ if and only if } k = \pm s.
\end{align*}
Moreover, the transversality condition is satisfied, that is,
\begin{align*}
	\partial_\lambda DF(\b{0}, 0)[\b{q}] \not\in \mathrm{ran} \,  DF(\b{0}, 0).
\end{align*}
A corresponding statement holds true for $DG(\b{0}, 0): \mathcal{X}_1 \to \mathcal{X}_1$.
\end{prop}

\subsection{The Proof of Theorem~\ref{thm_poly}}

Let us fix some $s \in \N$, and choose $(\tau_k)_{k \in \Z}$ as in~\eqref{eq_assumptions}. We introduce the trivial family $\mathcal{T} := \{ (\b{0}, \lambda) \in \mathcal{X}_1 \times \R
\: | \: \lambda \in \R \}$.

\begin{steps}
\step{Proof of (i).}
By Proposition~\ref{prop_mapping}, the maps $F$ resp. $G$ are smooth and vanish on the trivial family $\mathcal{T}$. In view of Proposition~\ref{prop_kernel}, the Crandall-Rabinowitz Theorem shows that $(\b{0}, 0) \in \mathcal{T}$ is a bifurcation point for $F(\b{v}, \lambda) = 0$ resp. $G(\b{v}, \lambda) = 0$ and provides an open interval $J_s \subseteq \R$ containing $0$ and a smooth curve 
\begin{align*}
	J_s \to \mathcal{X}_1 \times \R,
	\qquad
	\alpha \mapsto (\b{v}^\alpha, \lambda^\alpha) = \left( (v_k^\alpha)_{k \in \Z}, \lambda^\alpha \right)
\end{align*}
of zeros of $F$ resp. $G$
(we do not denote its dependence on $s$)
with $\b{v}^0 = \b{0}, \lambda^0 = 0$ as well as $\frac{\mathrm{d}}{\mathrm{d}\alpha} \big|_{\alpha = 0} \b{v}^\alpha = \b{q}$ where $\b{q}$ is a nontrivial element of the kernel of $DF(\b{0}, 0)$ resp. $DG(\b{0}, 0)$. 
We let $\b{u}^\alpha := \b{v}^\alpha + \b{w}$ and define polychromatic states $U^\alpha$ as in (i). Then $U^\alpha$ is a classical solution of the cubic Klein-Gordon equation~\eqref{eq_wave} due to Proposition~\ref{prop_mapping} since $F(\b{v}^\alpha, \lambda^\alpha) = 0$ resp. $G(\b{v}^\alpha, \lambda^\alpha) = 0$. By their very definition, the solutions $U^\alpha$ are time-periodic with period $2 \pi / \omega$ (maybe less). This proves (i).

\step{Proof of (ii).}
Since $F$ resp. $G$ are smooth, so is the map $J_s \to \mathcal{X}_1 \times \R, \: \alpha \mapsto (\b{v}^\alpha, \lambda^\alpha)$.
By Proposition~\ref{prop_kernel}, $q_k \neq 0$ if and only if $k = \pm s$, which implies that only the $\pm s$-th components of
\begin{align*}
	\frac{\mathrm{d}}{\mathrm{d}\alpha} \bigg|_{\alpha = 0} \b{u}^\alpha
	= \frac{\mathrm{d}}{\mathrm{d}\alpha} \bigg|_{\alpha = 0} \b{v}^\alpha = \b{q}
\end{align*}
do not vanish.
For sufficiently small nonzero values of $\alpha$, the solutions $U^\alpha$ are thus nonstationary.
In particular, the direction of bifurcation changes when changing the value of $s$, and the associated bifurcating curves are, at least locally, mutually different.  

\step{Proof of (iii).}
We show finally that, under the additional assumption that $\Gamma(x) \neq 0$ for almost all $x \in \R^3$, every non-stationary solution
\begin{align*}
	U^\alpha(t, x) = w_0(x) + v_0^\alpha(x) + \sum_{k=1}^\infty 2 \cos(\omega k t) \: v_k^\alpha(x)
\end{align*}
in fact possesses infinitely many nontrivial coefficients $v_k^\alpha$. 
Indeed, assuming the contrary, we can choose a maximal $r > 0$ (since $U^\alpha$ is non-stationary) with $v_r^\alpha \not\equiv 0$ or equivalently $u_r^\alpha = v_r^\alpha + w_r  = v_r^\alpha \not\equiv 0$. But then, 
\begin{align*}
	v_{3r}^\alpha 
	= 	\sum_{l+m+n = 3r} \mathcal{R}_{\mu_{3r}}^{\tau_{3r}} [\Gamma \: u_l^\alpha \, u_m^\alpha \, u_n^\alpha]
	= \mathcal{R}_{\mu_{3r}}^{\tau_{3r}} [\Gamma \: (v_r^\alpha)^3] \not\equiv 0
\end{align*}
since the convolution identity implies $- \Delta v_{3r}^\alpha - \mu_{3r} v_{3r}^\alpha = \Gamma \: (v_r^\alpha)^3$, and $\Gamma \: (v_r^\alpha)^3 \not\equiv 0$ since $\Gamma(x) \neq 0$ almost everywhere by assumption. This contradicts the maximality of $r$.
\hfill $\square$
\end{steps}

\subsection{The Proof of Remark~\ref{rmk} (c)}

Finally, as announced in Remark~\ref{rmk}~(c), we verify the nondegeneracy assumption~\eqref{eq_nondegenerate,2} resp.~\eqref{eq_nondegenerate,S} for constant positive $\Gamma$. 

\begin{lem}[Nondegeneracy, \`{a} la Bates and Shi~\cite{bates}]\label{lem_nondegeneracy}
Let $\Gamma \equiv \Gamma_0$ for some $\Gamma_0 > 0$, and assume that $w_0 \in C^2_\mathrm{rad}(\R^3)$ is a radially symmetric solution of~\eqref{eq_w0} the profile of which satisfies $w_0(r) > 0$, $w_0'(r) < 0$ for all $r>0$, and both $w_0(r)$ and $w_0'(r)$ decay exponentially as $r \to \infty$. Then the nondegeneracy property~\eqref{eq_nondegenerate,S} holds, i.e. for any radial, twice differentiable $q_0 \in X_1$
\begin{align*}
	- \Delta q_0 + q_0 = 3 \Gamma_0 \: w_0^2 \: q_0
	\text{ on } \R^3
	\qquad \text{implies } 
	\qquad
	q_0 \equiv 0.
\end{align*}
\end{lem}
This can be proved closely following the line of argumentation by Bates and Shi~\cite{bates}, Theorem~5.4~(6). The main difference is that they state the nondegeneracy result as a spectral property of the operator $- \Delta + m^2 + 3 \Gamma_0 w_0^2: \: H^2(\R^3) \to L^2(\R^3)$ whereas we cannot use the Hilbert space setting but discuss solutions in $X_1$. However, the technique of Bates and Shi (and also of Wei's proof in~\cite{wei}) is based on an expansion at a fixed radius $r > 0$ in terms of the eigenfunctions of the Laplace-Beltrami operator on $L^2(\mathbb{S}^2)$. This provides coefficients depending on $r$, and the conclusions are obtained from the analysis of these profiles on an ODE level using results due to Kwong and Zhang~\cite{kwong}. These ideas apply in the topology of $X_1$ in the very same way; for details, cf.~\cite{myDiss}, (proof of) Lemma~4.11.

\section{Proofs of the Auxiliary Results}~\label{ch_wave-proofs}

\vspace*{-0.5cm}

\begin{proofof}{Proposition~\ref{prop_convolution-x1}}
Let $\b{u}^{(j)} = (u^{(j)}_k)_{k \in \Z} \in \mathcal{X}_1$ for $j= 1, 2, 3$.  We find the following chain of inequalities
\begin{align*}
	&\norm{\b{u}^{(1)} \star \b{u}^{(2)} \star \b{u}^{(3)}}_{\mathcal{X}_3}
	=
	\sum_{k \in \Z} \norm{(\b{u}^{(1)} \star \b{u}^{(2)} \star \b{u}^{(3)})_k}_{X_3}
		\\
	&\quad \leq 
	\sum_{k \in \Z} \sum_{\substack{l, m, n \in \Z \\ l + m + n = k}} \norm{ 	u_l^{(1)} \, u_m^{(2)} \, u_n^{(3)} }_{X_3}
	\\
	&\quad \leq 
	\sum_{k \in \Z} \sum_{\substack{l, m, n \in \Z \\ l + m + n = k}} 
	\norm{u_l^{(1)}}_{X_1} \, \norm{u_m^{(2)}}_{X_1} \, \norm{u_n^{(3)}}_{X_1}
	\\
	&\quad= 
	\norm{\left( \norm{u_l^{(1)}}_{X_1} \right)_{l \in \Z} \star \left( \norm{u_m^{(2)}}_{X_1} \right)_{m \in \Z}
	\star \left( \norm{u_n^{(3)}}_{X_1} \right)_{n \in \Z}}_{\ell^1(\Z)}
	\\
	&\quad\leq 
	\norm{\left( \norm{u_l^{(1)}}_{X_1} \right)_{l \in \Z}}_{\ell^1(\Z)} 
	\norm{\left( \norm{u_m^{(2)}}_{X_1} \right)_{m \in \Z}}_{\ell^1(\Z)} 
	\norm{\left( \norm{u_n^{(3)}}_{X_1} \right)_{n \in \Z}}_{\ell^1(\Z)} 
	\\
	&\quad = 
	\norm{\b{u}^{(1)}}_{\mathcal{X}_1} \norm{\b{u}^{(2)}}_{\mathcal{X}_1} \norm{\b{u}^{(3)}}_{\mathcal{X}_1},
\end{align*}
where finally Young's inequality for convolutions in $\ell^1(\Z)$ has been applied. 
Since the latter term is finite, we infer
$\b{u}^{(1)} \star \b{u}^{(2)} \star \b{u}^{(3)} \in \mathcal{X}_3$.
\end{proofof}

\begin{proofof}{Proposition~\ref{prop_mapping}}
\begin{steps}
\step{Decay estimates}

The proof of Proposition~\ref{prop_mapping} requires convergence properties in order to handle the infinite series in the definition of $U(t,x)$, which we first provide in the following two lemmas.

\begin{lem}\label{lem_scaling}
The convolution operators $\mathcal{R}_\mu^\tau: X_3 \to X_1$ satisfy for $\tau \in (0, \pi)$ and $\mu > 0$
\begin{equation*}
	\qquad \forall \: f \in X_3
	\!\!\!\!\!\!\!\!\!\!\!\!\!\!\!\!\!\!\!\!\!\!\!\!\!\!\!\!\!\!\!\!\!\!\!\!\!\!\!\!
	\begin{split}
	&\norm{\mathcal{R}_\mu^\tau [f]}_{X_1} 
	\leq \frac{C}{\sin(\tau)} \: \left( 1+\frac{1}{\sqrt{\mu}} \right) \:  \cdot \norm{f}_{X_3},
	\\
	&\norm{\mathcal{R}_\mu^\tau [f]}_{L^4(\R^3)} 
	\leq \frac{C}{\sqrt[4]{\mu} \: \sin(\tau)} \cdot \norm{f}_{L^\frac{4}{3}(\R^3)}.
	\end{split}
\end{equation*}
\end{lem} 

The fact that a power of $\mu$ appears in the denominator is crucial since it will finally provide the convergence and regularity of the polychromatic sums where $\mu = \mu_k = \omega^2 k^2 - m^2$ for $k \in \Z$.
\\
The proof of Lemma~\ref{lem_scaling} relies, via rescaling, on the respective estimates for $\mu = 1$. These can be found in~\cite{own_cubic}, pp.~1038--1039 for the $X_3$-$X_1$ estimate and in~\cite{EvequozWeth},~Theorem~2.1 for the $L^{4/3}$-$L^4$ estimate.

\begin{lem}\label{lem_uniformdecay}
Let $\Gamma \in L^\infty_\mathrm{rad}(\R^3) \cap C^1_\mathrm{loc}(\R^3)$ and assume $\b{u} = (u_k)_{k\in\Z} \in \mathcal{X}_1$ is a sequence of $C^2_\mathrm{loc}$ functions which satisfy the following system of convolution equations:
\begin{align*}
	u_k = \mathcal{R}_{\mu_k}^{\tau_k}  [\Gamma \: (\b{u} \star \b{u} \star \b{u})_k]
	\qquad \qquad \qquad \text{for all } k \in \Z \text{ with } |k| > s
\end{align*}
where $\mu_k = \omega^2 k^2 - m^2$ and $\tau_k \in (\delta, \pi - \delta)$ for some $\omega > m, \delta \in \left(0, \frac{\pi}{2} \right)$. Then there holds:
\begin{itemize}
\item[(i)]
For every $\alpha \geq 0$, there exists a constant $C_\alpha \geq 0$ with
\begin{align*}
	\norm{u_k}_{L^4(\R^3)} + \norm{\Gamma \: (\b{u} \star \b{u} \star \b{u})_k}_{L^4(\R^3)} 
	\leq C_\alpha \cdot (k^2 + 1)^{-\frac{\alpha}{2}}
	\qquad
	(k \in \Z).
\end{align*} 
\item[(ii)]
For every ball $B = B_R(0) \subseteq \R^3$ and $\alpha \geq 0$ there exists a constant $D_\alpha(B) \geq 0$ with
\begin{align*}
	|u_k(x)| + |\nabla u_k(x)| + |D^2 u_k(x)| \leq D_\alpha(B) \cdot (k^2 + 1)^{-\frac{\alpha}{2}}
	\qquad
	(k \in \Z, x \in B).
\end{align*}
\item[(iii)]
For every $\alpha \geq 0$, there exists a constant $E_\alpha \geq 0$ with
\begin{align*}
	\norm{u_k}_{X_1} \leq E_\alpha \cdot (k^2 + 1)^{-\frac{\alpha}{2}}
	\qquad
	(k \in \Z).
\end{align*} 
\end{itemize}
\end{lem}

The proof of Lemma~\ref{lem_uniformdecay} can also be found in detail in the author's PhD thesis~\cite{myDiss}, Lemma~4.13. We present here the most important step and summarize the remainder briefly, since it is mainly based on the application of (standard) elliptic regularity estimates. 

As $u_k \in X_1 \cap C^2_\text{loc}(\R^3)$ for all $k \in \Z$ by assumption, it is straightforward to find constants as in the lemma for a finite number of elements $u_{-s}, ..., u_s$.  Hence it is sufficient to study those $k \in \Z$ with $|k| > s$; for these, we have $\mu_k = k^2 \omega^2 - m^2 \geq c_s (k^2 +1)$ for some positive $c_s > 0$ depending on the parameters $\omega$ and $m$.
The decay estimates of arbitrary order in $k$ we aim to prove essentially go back to the $L^{4/3}$-$L^4$ scaling property stated in Lemma~\ref{lem_scaling} above. Indeed, due to $\delta < \tau_k < \pi - \delta$, it provides $C_1 = C_1(\norm{\Gamma}_\infty, \delta, \omega, m, s) \geq 0$ with
\begin{align}\label{eq_proof-scale-k}
	\norm{u_k}_{L^4(\R^3)} 
	\leq 
	\frac{C_1}{(k^2 + 1)^\frac{1}{4}} \: 
	\norm{(\b{u} \star \b{u} \star \b{u})_k}_{L^\frac{4}{3}(\R^3)}
	\qquad \text{for all } k \in \Z.
\end{align}
With that, assuming $\sum_{k \in \Z} (k^2 + 1)^\frac{\alpha}{2} \norm{u_k}_{L^4(\R^3)} < \infty$ for some $\alpha \geq 0$ (which is trivially satisfied for $\alpha = 0$ since $\b{u} \in \mathcal{X}_1$), one can iterate as follows
\begin{align*}
	&\sum_{k \in \Z} (k^2 + 1)^\frac{\alpha + 1/2}{2} \: \norm{u_k}_{L^4(\R^3)}
	\overset{\eqref{eq_proof-scale-k}}{\leq} 
	C_1 \: \sum_{k \in \Z} (k^2 + 1)^\frac{\alpha}{2}
	\: \norm{(\b{u} \star \b{u} \star \b{u})_k}_{L^\frac{4}{3}(\R^3)}
	 \\
	& \quad \leq C_1 \:  
	 \sum_{k \in \Z} \sum_{l+m+n = k} 
	 ((l+m+n)^2 + 1)^\frac{\alpha}{2} \norm{u_l}_{L^4(\R^3)} \, \norm{u_m}_{L^4(\R^3)} \, \norm{u_n}_{L^4(\R^3)} 
	 \\
	&  \quad \leq 2^\alpha \: C_1 \:  
	 \sum_{k \in \Z} \sum_{l+m+n = k} 
	 \left[ (l^2 + 1)^\frac{\alpha}{2} \norm{u_l}_{L^4(\R^3)} \, (m^2 + 1)^\frac{\alpha}{2} \norm{u_m}_{L^4(\R^3)} \, (n^2 + 1)^\frac{\alpha}{2} \norm{u_n}_{L^4(\R^3)} \right]
	\\
	& \quad = 2^\alpha \: C_1 \:  
	 \left( \sum_{k \in \Z} (k^2 + 1)^\frac{\alpha}{2} \norm{u_k}_{L^4(\R^3)} \right)^3
	 \\
	 & \quad < \infty.
\end{align*}
This shows the first part of the estimate in (i), and the second part follows by combining the former with the interpolation estimate
\begin{align*}
	\norm{u_l u_m u_n}_{L^4(\R^3)}
	&\leq \norm{u_l}_{L^{12}(\R^3)}\norm{u_m}_{L^{12}(\R^3)}\norm{u_n}_{L^{12}(\R^3)}
	\\
	&\leq \left[ \norm{u_l}_{L^{4}(\R^3)}\norm{u_m}_{L^{4}(\R^3)}\norm{u_n}_{L^{4}(\R^3)} \right]^\frac{1}{3} \: \Big[ \norm{u_l}_{\infty}\norm{u_m}_{\infty}\norm{u_n}_{\infty} \Big]^\frac{2}{3}
	\\
	&\leq \left[ \norm{u_l}_{L^{4}(\R^3)}\norm{u_m}_{L^{4}(\R^3)}\norm{u_n}_{L^{4}(\R^3)} \right]^\frac{1}{3} \: \norm{\b{u}}^2_{\mathcal{X}_1}.
\end{align*}
The local estimate in (ii) can be derived from the global $L^4$ bounds in (i) using elliptic regularity, which first provides estimates in $W^{2,4}_\text{loc}(\R^3)$ and then is suitable Hölder spaces. The estimate (iii) in the $X_1$ norm essentially uses the explicit representations (given $f \in X_3$)
\begin{align*}
	&\mathcal{R}_\mu^\tau [f](x) 
	= \int_{\R^3}  \frac{\sin(|x-y|\sqrt{\mu_k} + \tau_k)}{4 \pi |x-y| \sin(\tau_k)} 
	\cdot f(y)  
	\: \mathrm{d}y
	\\
	& = \frac{\sin(|x| \sqrt{\mu_k} + \tau_k)}{|x| \sin(\tau_k)}   \int_0^{|x|}	\frac{\sin(r \sqrt{\mu_k}) }{r \sqrt{\mu_k}} 	
 	 f(r) \: r^2 \:  \mathrm{d}r
	+ \frac{\sin(|x| \sqrt{\mu_k})}{|x| \sin(\tau_k)} 
	 \int_{|x|}^\infty \frac{\sin(r \sqrt{\mu_k} + \tau_k) }{r \sqrt{\mu_k}}  f(r) \: r^2 	\:  \mathrm{d}r.
\end{align*}
Starting here, Hölder's inequality and (i) yield (iii); again, for details, cf.~\cite{myDiss}.

\step{Mapping properties of $F$ resp. $G$.}
For $\lambda \in \R$ and $\b{v} \in \mathcal{X}_1$, we set $\b{u} := \b{w} + \b{v}$ and recall the defining equations~\eqref{eq_F}~and~\eqref{eq_G}:
\begin{align*}
	&F(\b{v}, \lambda)_{k}
:= v_k -
\begin{cases}
	\mathcal{P}_{\mu_0}
	\left[ \Gamma \: \left(\b{u} \star \b{u} \star \b{u}\right)_0 
	- \Gamma \: w_0^3  \right] 
	& k = 0,
	\\
	\mathcal{R}_{\mu_s}^{\pi/2}
	\left[ \Gamma \: \left(\b{u} \star \b{u} \star \b{u}\right)_{\pm s} \right]
	\\
	\qquad 
	+ (\cot(\tau_{\pm s}) - \lambda) \tilde{\Psi}_{\mu_s} \ast
	\left[ \Gamma \: \left(\b{u} \star \b{u} \star \b{u}\right)_{\pm s} \right]
	& k = \pm s,
	\\
	\mathcal{R}_{\mu_k}^{\tau_{k}}
	\left[ \Gamma \: \left(\b{u} \star \b{u} \star \b{u}\right)_{k} \right]
	& \text{else};
\end{cases}
\\
	&G(\b{v}, \lambda)_{k}
:= v_k - 
\begin{cases}
	\mathcal{P}_{\mu_0}
	\left[ \Gamma \: \left(\b{u} \star \b{u} \star \b{u}\right)_0 
	-  \Gamma \: w_0^3  \right] 
	& k = 0,
	\\
	\mathcal{R}_{\mu_s}^{\pi/2}
	\left[ \Gamma \: \left(\b{u} \star \b{u} \star \b{u}\right)_{\pm s} \right] 
	\\
	\qquad
	+ (1 - \lambda)
	\left( \alpha^{(\mu_s)}(v_{\pm s}) + \beta^{(\mu_s)}(v_{\pm s}) \right) \tilde{\Psi}_{\mu_s}
	& k = \pm s, 
	\\
	\mathcal{R}_{\mu_k}^{\tau_{k}}
	\left[ \Gamma \: \left(\b{u} \star \b{u} \star \b{u}\right)_{k} \right]
	& \text{else}.
\end{cases}
\end{align*}
Our main concern will be convergence of the infinite sums related to the space $\mathcal{X}_1 = \ell^1_\mathrm{sym}(\Z, X_1)$. Noticing that $F$ and $G$ only differ in the $\pm s$-th component, and that the scalar parameter $\lambda$ only appears as a multiplicative factor, we solely discuss smoothness of the map $F(\,\cdot\, , \lambda): \mathcal{X}_1 \to \mathcal{X}_1$ with $\lambda \in \R$ fixed.

The main tool is the following uniform norm estimate for the operators appearing in the components of $F$.
Recalling that $\tau_k \in \{ \frac{\pi}{4}, \frac{3\pi}{4} \}$ for $k \neq 0, \pm s$ by~\eqref{eq_assumptions}, Lemma~\ref{lem_scaling} above (for $k \neq 0, \pm s$) as well as the continuity properties stated in Lemmas~\ref{lem_convol} and~\ref{lem_linS} (for $k = \pm s$ and $k = 0$, respectively) provide a constant $C_0 = C_0(\lambda, \tau_s, \omega, m) > 0$ with
\begin{equation}\label{eq_proof-C0}
\begin{split}
	&\norm{\mathcal{R}_{\mu_k}^{\tau_{k}}}_{\mathcal{L}(X_3, X_1)} \leq C_0 \quad (k \in \Z \setminus \{ \pm s \}), 
	\\
	&\norm{\mathcal{R}_{\mu_s}^{\pi/2}}_{\mathcal{L}(X_3, X_1)} \leq \frac{C_0}{2},
	\:\:
	\norm{(\cot(\tau_{\pm s}) - \lambda) \: \tilde{\Psi}_{\mu_s} \ast}_{\mathcal{L}(X_3, X_1)} \leq \frac{C_0}{2},
	\\
	&\norm{\mathcal{P}_{\mu_0}}_{\mathcal{L}(X_3, X_1)} \leq C_0.
\end{split}
\end{equation}

We now let $\b{v} \in \mathcal{X}_1$ and define $\b{u} = \b{v} + \b{w}$. 
Since $\Gamma$ is assumed to be continuous and bounded, Proposition~\ref{prop_convolution-x1} implies that $\Gamma \: \left(\b{u}\star\b{u}\star\b{u}\right) \in \mathcal{X}_3$. Thus every component $F(\b{v}, \lambda)_k$ is a well-defined element of $X_1$, and we estimate
\begin{align*}
	\norm{F(\b{v}, \lambda)}_{\mathcal{X}_1}
	&= \sum_{k \in \Z} \norm{F(\b{v}, \lambda)_k}_{X_1}
	\\
	&  \overset{\eqref{eq_proof-C0}}{\leq}
	\norm{\b{v}}_{\mathcal{X}_1} 
	+ C_0 \norm{\Gamma w_0^3}_{X_3}
	 + C_0 \sum_{k \in \Z} \norm{\Gamma \left(\b{u} \star \b{u} \star \b{u}\right)_k}_{X_3}
	 \\
	 & \!\!\!   \overset{\text{Prop.}~\ref{prop_convolution-x1}}{\leq}
	\norm{\b{v}}_{\mathcal{X}_1} + 
	C_0 \norm{\Gamma}_\infty \norm{w_0}_{X_1}^3 
	+ C_0 \norm{\Gamma}_\infty \norm{\b{u}}_{\mathcal{X}_1}^3.
\end{align*}
This is finite, hence $F(\b{v}, \lambda) \in \mathcal{X}_1$ as asserted.
Since $F(\,\cdot\, , \lambda)$ is a combination of continuous linear operators and polynomials in the convolution algebra, essentially the same estimates can be used to show differentiability (to arbitrary order); one thus obtains in particular~\eqref{eq_DF}.

\step{Solution properties of $u_k(x)$.}
First of all, recalling that $\b{w} = (..., 0, w_0, 0, ...)$ and hence $(\b{w} \star  \b{w} \star \b{w})_k = \delta_{k, 0} \, w_0^3$ for $k \in \Z$, one can immediately see that $F(\b{0},\lambda) = G(\b{0},\lambda) = \b{0}$ for all $\lambda \in \R$.
Let us now assume that $F(\b{v}, \lambda) = 0$ resp. $G(\b{v}, \lambda) = 0$ for some $\b{v} \in \mathcal{X}_1$ and $\lambda \in \R$. 
Again, we define $\b{u} := \b{v} + \b{w}$, and summarize
\begin{align*}
	u_0 - w_0 &= v_0 
	=	\mathcal{P}_{\mu_0} \left[ \Gamma  \left(\b{u} \star \b{u} \star \b{u}\right)_0 
	- \Gamma \: w_0^3  \right],
	\\
	u_{\pm s} &= v_{\pm s} 
	= \mathcal{R}_{\mu_s}^{\pi/2}
	\left[ \Gamma  \left(\b{u} \star \b{u} \star \b{u}\right)_{\pm s} \right]
	+ \begin{cases}
	(\cot(\tau_s) - \lambda)  \tilde{\Psi}_{\mu_s} \ast
	\left[ \Gamma  \left(\b{u} \star \b{u} \star \b{u}\right)_{\pm s} \right],
	\\
	 (1 - \lambda) 
	\left( \alpha^{(\mu_s)}(v_{\pm s}) + \beta^{(\mu_s)}(v_{\pm s}) \right) \tilde{\Psi}_{\mu_s},
	\end{cases}
	\\
	u_k &= v_k 
	= \mathcal{R}_{\mu_k}^{\tau_k}
	\left[ \Gamma  \left(\b{u} \star \b{u} \star \b{u}\right)_k \right]
	\qquad (k \in \Z \setminus \{ 0, \pm s \}).
\end{align*}
By choice of $\tau_k$ in equation~\eqref{eq_assumptions}, we observe in particular that the requirements of Lemma~\ref{lem_uniformdecay} are satisfied with any $\delta < \frac{\pi}{4}$, which we will rely on throughout the subsequent steps. 
But first, according to Lemmas~\ref{lem_convol}~and~\ref{lem_linS}, $v_k, u_k \in X_1 \cap C^2_\mathrm{loc}(\R^3)$ satisfy the differential equations
\begin{align*}
	- \Delta v_k - \mu_k v_k = \Gamma (x) \: \big[ (\b{u}\star\b{u}\star\b{u})_k - \delta_{k,0} w_0^3 \big]
	\quad \text{on } \R^3
\end{align*}
or equivalently, in view of $\b{w} = (..., 0, w_0, 0, ...)$, of~\eqref{eq_w0} and of $\mu_k = \omega^2 k^2 - m^2$,
\begin{align}\label{eq_stationary-proof}
	- \Delta u_k - (\omega^2 k^2 - m^2) u_k = \Gamma (x) \: (\b{u}\star\b{u}\star\b{u})_k 
	\quad \text{on } \R^3.
\end{align}
We now define formally for $t \in \R, x \in \R^3$ 
\begin{align}\label{eq_poly-proof}
	U(t, x) := w_0(x) + v_0(x) + \sum_{k=1}^\infty 2 \cos(\omega kt) \: v_k(x)
	= \sum_{k \in \Z} \mathrm{e}^{\i \omega k t} \: u_k(x).
\end{align}
Since by assumption $\b{u} = \b{v} + \b{w} \in \ell^1(\Z, X_1)$, the Weierstrass M-test asserts that the sum in~\eqref{eq_poly-proof} converges in $X_1$ uniformly with respect to $t \in \R$, and hence the map $t \mapsto U(t, \,\cdot \,)$ is continuous as a map from $\R$ to $X_1$. We next show stronger regularity properties of $U(t,x)$.

\step{Differentiability of $U(t, x)$.}
We prove that the map $t \mapsto U(t, \,\cdot \,)$, when interpreted as a map from $\R$ to $X_1$, possesses two continuous time derivatives given by
\begin{align*}
	\partial_t U(t, \,\cdot \,) = \sum_{k \in \Z} \i \omega k \: \mathrm{e}^{\i \omega k t} \: u_k,
	\qquad
	\partial_t^2 U(t, \,\cdot \,) = \sum_{k \in \Z} -\omega^2 k^2 \: \mathrm{e}^{\i \omega k t} \: u_k.
\end{align*}
Indeed, term-by-term differentiation is justified since the sums above as well as in~\eqref{eq_poly-proof} converge in $X_1$ uniformly with respect to time. This is a consequence of the Weierstraß M-test and the decay estimate in Lemma~\ref{lem_uniformdecay}~(iii).
Hence, as asserted, the map $t \mapsto U(t, \,\cdot \,)$ is twice continuously differentiable as a map from $\R$ to $X_1$ - the same strategy yields in fact $C^\infty$ regularity in time.

Similarly, the local regularity estimate in Lemma~\ref{lem_uniformdecay}~(ii) implies $U \in C^2(\R \times B)$ for every given ball $B = B_R(0) \subseteq \R^3$ again via term-by-term differentiation.   
Since the radius of the ball $B$ is arbitrary, we conclude for $t \in \R$ and all $x \in \R^3$
\begin{align*}
	\left[ \partial_t^2 - \Delta + m^2 \right] U(t, x)
	&=
	\sum_{k \in \Z}  \mathrm{e}^{\i \omega k t} \: \left[ -\omega^2 k^2 - \Delta + m^2 \right] u_k(x) 
	\\
	&\!\!\overset{\eqref{eq_stationary-proof}}{=}
	\sum_{k \in \Z}  \mathrm{e}^{\i \omega k t} \: \Gamma(x) \: \sum_{l+m+n = k} u_l(x) \, u_m(x) \, u_n(x)
	\\
	 &=
	\Gamma(x) \: \left( \sum_{l \in \Z}  \mathrm{e}^{\i \omega l t} \: u_l(x) \right)
	\left( \sum_{m \in \Z}  \mathrm{e}^{\i \omega m t} \: u_m(x) \right)
	\left( \sum_{n \in \Z}  \mathrm{e}^{\i \omega n t} \: u_n(x) \right)
	\\
	&= \Gamma(x) \: U(t,x)^3
\end{align*}
where the re-ordering of the summation is justified by absolute convergence of the sums.
Thus $U$ is shown to be a classical solution of the Klein-Gordon equation~\eqref{eq_wave}. 
\end{steps}
\end{proofof}

\begin{proofof}{Proposition~\ref{prop_kernel}}
We prove the statement for the map $F$ and then comment on the aspects that differ in case of $G$.
Using formula~\eqref{eq_DF}, we find for $k \in \Z$ and $\b{q} \in \mathcal{X}_1$, recalling that $w_k = 0$ for $k \in \Z \setminus \{ 0 \}$ and that $\mathcal{R}_{\mu_s}^{\tau_s} = \mathcal{R}_{\mu_s}^{\pi/2} + \cot(\tau_s) \, \tilde{\Psi}_{\mu_s} \ast$,
\begin{align*}
	&DF(\b{0}, 0)[\b{q}]_k = q_k - 3 \: \mathcal{R}_{\mu_k}^{\tau_k} [\Gamma \: (\b{q} \star \b{w} \star \b{w})_k]
	= q_k - 3 \: \mathcal{R}_{\mu_k}^{\tau_k} \left[ \Gamma \: w_0^2 \cdot q_k \right],
	\\
	&DF(\b{0}, 0)[\b{q}]_0 = q_0 - 3 \: \mathcal{P}_{\mu_0} [\Gamma \: (\b{q} \star \b{w} \star \b{w})_0] = q_0 - 3 \: \mathcal{P}_{\mu_0} \left[ \Gamma \: w_0^2 \cdot q_0 \right].
\end{align*}
For $\b{q} \in \ker DF(\b{0}, 0)$, and in view of the choice of $\tau_k$ in~\eqref{eq_assumptions}, the nondegeneracy properties~\eqref{eq_nondegenerate} imply $q_k \equiv 0$ for $k \in \Z, k \neq \pm s$. 
Since $\tau_{\pm s} = \sigma_{s}$ in~\eqref{eq_assumptions}, Proposition~\ref{prop_chapter2} guarantees the existence of a nontrivial solution $q_s \in X_1$ of
\begin{align}\label{eq_qs-case1}
	q_s = 3 \: \mathcal{R}_{\mu_s}^{\tau_s} \left[ \Gamma \: w_0^2 \cdot q_s \right]
\end{align}
which is unique up to a multiplicative factor. Hence $\ker DF(\b{0}, 0)$ has the asserted form.  (We recall here that we consider the subspace of symmetric sequences, whence $q_{-s} = q_s$.)
Further, by Lemmas~\ref{lem_convol}~and~\ref{lem_linS} in the final Section~\ref{sect_helmholtz}, the operators
\begin{align*}
	X_1 \to X_1, \qquad
	\begin{cases}
	q_k \mapsto q_k - 3 \: \mathcal{R}_{\mu_k}^{\tau_k} \left[ \Gamma \: w_0^2 \cdot q_k \right]
	& (k \neq 0)
	\\
	q_0 \mapsto q_0 - 3 \: \mathcal{P}_{\mu_0} \left[ \Gamma \: w_0^2 \cdot q_0 \right]
	\end{cases}
\end{align*}
are linear compact perturbations of the identity and so $\ker DF(\b{0}, 0)$ is 1-1-Fredholm. In order to verify transversality, we compute for $k \in \Z$ and $\b{q} \in \ker DF(\b{0}, 0) \setminus \{ \b{0} \}$
\begin{align*}
	\partial_\lambda DF(\b{0}, 0)[\b{q}]_k
	= \begin{cases}
		3 \: \tilde{\Psi}_{\mu_s} \ast [ \Gamma \: w_0^2 \, q_s ],		& k = \pm s,
		\\
		0, & \text{else}.
	\end{cases} 
\end{align*}
Assuming for contradiction that $\partial_\lambda DF(\b{0}, 0)[\b{q}] = DF(\b{0}, 0)[\b{p}]$ for some $\b{p} \in \mathcal{X}_1$, we infer in particular that the component $p_s$ satisfies the convolution identity
\begin{align}\label{eq_aux-transversal}
	p_s - 3 \: \mathcal{R}_{\mu_s}^{\tau_s} \left[ \Gamma \: w_0^2 \cdot p_s \right] 
	=  3  \: \tilde{\Psi}_{\mu_s} \ast  [ \Gamma \: w_0^2 \cdot q_s ]
\end{align}
and hence, following Lemmas~\ref{lem_convol},~\ref{lem_linH}
\begin{align*}
	- \Delta p_s - \mu_s p_s = 3 \: \Gamma(x) \:  w_0^2(x) \, p_s
	\quad \text{on } \R^3,
\end{align*}
which is also nontrivially solved by $q_s$ as a consequence of~\eqref{eq_qs-case1}. Due to the uniqueness statement in Proposition~\ref{prop_chapter2}, this implies that $p_s = c \cdot q_s$ for some $c \in \R$. But then, applying~\eqref{eq_qs-case1} to~\eqref{eq_aux-transversal}, we obtain 
$	\tilde{\Psi}_{\mu_s} \ast [ \Gamma \: w_0^2 \cdot q_s ] = 0$. Hence by the asymptotic expansion in Lemma~\ref{lem_convol} 
\begin{align*}
	\widehat{\Gamma w_0^2 q_s} (\sqrt{\mu_s}) = 0
\end{align*}
and therefore, due to $q_s = 3 \: \mathcal{R}^{\tau_s}_{\mu_s} [ \Gamma \: w_0^2 \, q_s ]$ and  Lemma~\ref{lem_linH},
\begin{align*}
	q_s(x) = O\left(\frac{1}{|x|^2}\right) \text{  as } |x| \to \infty.
\end{align*}
This contradicts Proposition~\ref{prop_chapter2} stating that the leading-order term as $|x| \to \infty$ of a nontrivial solution $q_s$ of $- \Delta q_s - \mu_s q_s = 3 \: \Gamma(x) \:  w_0^2(x) \, q_s$ cannot vanish.

In the case $\tau_s = 0$, we see as above that $\b{q} \in \ker DG(\b{0}, 0)$ if and only if $q_k = 0$ for $k \neq \pm s$, and that $q_{s} = q_{-s}$ can be chosen to be the (nontrivial) solution of 
\begin{align}\label{eq_qs-case2}
	q_s = 3 \: \mathcal{R}_{\mu_s}^{\pi/2} \left[ \Gamma \: w_0^2 \cdot q_s \right] 
	+  \alpha^{(\mu_s)}(q_s)  \: \tilde{\Psi}_{\mu_s}
	\quad \text{with } \quad \beta^{(\mu_s)}(q_s) = 0.
\end{align}
Similarly, $\ker DG(\b{0}, 0)$ is 1-1-Fredholm. 
We again assume for contradiction that there is $\b{p} \in \mathcal{X}_1$ with $\partial_\lambda DG(\b{0}, 0)[\b{q}] = DG(\b{0}, 0)[\b{p}]$, which implies in particular
\begin{equation}\label{eq_aux-transversal-2}
	p_s - 3 \: \mathcal{R}_{\mu_s}^{\pi/2} \left[ \Gamma \: w_0^2 \cdot p_s \right] 
	- \left( \alpha^{(\mu_s)}(p_s) + \beta^{(\mu_s)}(p_s) \right)  \tilde{\Psi}_{\mu_s}
	=   \alpha^{(\mu_s)}(q_s)   \tilde{\Psi}_{\mu_s}
\end{equation}
with $\beta^{(\mu_s)}(q_s) = 0$. Thus, according to Lemma~\ref{lem_convol}, $p_s$ solves the differential equation
\begin{align*}
	- \Delta p_s - \mu_s p_s = 3 \: \Gamma(x) \:  w_0^2(x) \, p_s
	\quad \text{on } \R^3,
\end{align*}
which is also solved by $q_s$, see equation~\eqref{eq_qs-case2}. As before, the uniqueness property in Proposition~\ref{prop_chapter2} implies $p_s = c \cdot q_s$ for some $c \in \R$, and inserting this into the identity~\eqref{eq_aux-transversal-2}, comparison with~\eqref{eq_qs-case2} yields $\alpha^{(\mu_s)}(q_s) = 0$. Since also $\beta^{(\mu_s)}(q_s) = 0$, we infer from the definition of the functionals $\alpha^{(\mu_s)}, \beta^{(\mu_s)}$ preceding Lemma~\ref{lem_linH-0} that, again, $q_s(x) = O(1/|x|^2)$, contradicting Proposition~\ref{prop_chapter2}.
\end{proofof}

\section{Appendix: Stationary Linear Helmholtz and Schrödinger Equations}\label{sect_helmholtz}

Given $\mu > 0$, we study aspects of the solution theory of the linear equations
\begin{align}\label{eq_linSH}
	- \Delta u \pm \mu u = f			\qquad \text{on } \R^3.
\end{align}
In the case of a ``$+$'', equation~\eqref{eq_linSH} is said to be a Schrödinger equation. Given any right-hand side $f \in L^2(\R^3)$, a unique solution $u \in H^2(\R^3)$ can be obtained by applying the resolvent $(- \Delta + \mu)^{-1}$, which can be calculated explicitly by applying the Fourier transform
\begin{align*}
	u = (- \Delta + \mu)^{-1} f = \int_{\R^3} \frac{\hat{f}(\xi)}{|\xi|^2 + \mu} \: \mathrm{e}^{\mathrm{i} \skp{\,\cdot\,}{\xi}} \: \frac{\mathrm{d}\xi}{(2\pi)^{3/2}}.
\end{align*}
In the case of a Helmholtz equation, i.e. of a ``$-$'' sign in~\eqref{eq_linSH}, this is not possible since $\mu > 0$ belongs to the essential spectrum of $- \Delta$ on $\R^3$. A well-established strategy to find solutions in spaces other than $L^2(\R^3)$ is known as Limiting Absorption Principle(s). The idea is to replace $\mu$ by $\mu + \i \varepsilon$, apply an $L^2$-resolvent, and pass to the limit $\varepsilon \to 0$ in a suitable topology, i.e. formally
 \begin{align*}
	u = `` \lim_{\varepsilon \searrow 0}" \: (- \Delta - (\mu + \i \varepsilon))^{-1} f =
	`` \lim_{\varepsilon \searrow 0}" \:  \int_{\R^3} \frac{\hat{f}(\xi)}{|\xi|^2 - (\mu + \i \varepsilon)} \:  \mathrm{e}^{\mathrm{i} \skp{\,\cdot\,}{\xi}} \: \frac{\mathrm{d}\xi}{(2\pi)^{3/2}}.
\end{align*}
Using tools from harmonic analysis, such a construction of solutions of linear inhomogeneous Helmholtz equations has been successfully done by Agmon~\cite{agmon-scattering} in weighted $L^2$ spaces, and by Kenig, Ruiz and Sogge~\cite{Kenig} as well as Guti\'{e}rrez~\cite{Gutierrez} in certain pairs of $L^p$ spaces. The resolvent-type operator is, then, for sufficiently nice $f$, given by a convolution
\begin{align*}
	u = \frac{\mathrm{e}^{\i | \, \cdot \, | \sqrt{\mu}}}{4 \pi |\, \cdot \, |} \ast f.
\end{align*}
Such studies are completed by characterizations of the so-called Herglotz waves, i.e. the solutions of the homogeneous equation $- \Delta u - \mu u = 0$ on the respective spaces, see e.g.~\cite{agmon}.

\medskip

We study the case of (real-valued, radial) functions $f \in X_3, u \in X_1$ with the Banach spaces 
\begin{align*}
	X_q &:=
	\left\{  
	v \in C_\mathrm{rad}(\R^3) \big| \: \norm{v}_{X_q} := \norm{(1 + |\cdot|^2)^\frac{q}{2} v}_\infty < \infty 
	\right\},
	\qquad q \in \{ 1, 3 \}.
\end{align*}
These have been successfully applied in solving systems of cubic Helmholtz equations in~\cite{own_cubic}. Let us again point out that the decay rate prescribed in $X_1$ is the natural one for solutions of Helmholtz equations on the full space $\R^3$.
Such solutions of the Helmholtz equation 
\begin{align}\label{eq_linH}
	- \Delta u - \mu u = f			\qquad \text{on } \R^3
\end{align}
can be obtained using convolution operators with kernels $\Psi_\mu, \tilde{\Psi}_\mu$ given by
\begin{align*}
	\Psi_\mu(x) = \frac{\cos(|x|\sqrt{\mu})}{4 \pi |x|},
	\qquad
	\tilde{\Psi}_\mu(x) = \frac{\sin(|x|\sqrt{\mu})}{4 \pi |x|}
	\qquad
	(x \in \R^3 \setminus \{ 0 \}).
\end{align*}
Here $\Psi_\mu, \tilde{\Psi}_\mu$ are radial solutions of the homogeneous Helmholtz equation on $\R^3 \setminus \{ 0 \}$. We notice that $\tilde{\Psi}_\mu$ extends to a smooth solution of $- \Delta u - \mu u = 0$ in $X_1$ and it is, up to constant multiples, the only one. Moreover, the following holds:

\begin{lem}[\cite{own_cubic}, Proposition~4]\label{lem_convol}
The convolution operators $f \mapsto \Psi_\mu \ast f$, $f \mapsto \tilde{\Psi}_\mu \ast f$ are well-defined, linear and compact as operators from $X_3$ to $X_1$. Moreover, given $f \in X_3$, the functions $w := \Psi_\mu \ast f$ and $\tilde{w} := \tilde{\Psi}_\mu \ast f$ belong to $X_1 \cap C^2_{\mathrm{loc}}(\R^3)$ and satisfy
\begin{align*}
	&- \Delta w - \mu w = f \quad \text{on } \R^3,
	\qquad w(x) = 4\pi \: \sqrt{\frac{\pi}{2}} \: \hat{f}(\sqrt{\mu}) \: \: \Psi_\mu(x) + O\left( \frac{1}{|x|^2} \right);
	\\
	&- \Delta \tilde{w} - \mu \tilde{w} = 0 \quad \text{on } \R^3,
	\qquad \tilde{w}(x) = 4\pi \: \sqrt{\frac{\pi}{2}} \: \hat{f}(\sqrt{\mu}) \: \: \tilde{\Psi}_\mu(x).
\end{align*}
\end{lem}
Here $\hat{f}(\sqrt{\mu})$ refers to the profile of the Fourier transform on $\R^3$. Working in a radial setting with strongly decaying inhomogeneities $f \in X_3$, the properties in the previous Lemma (and in the following ones) can be verified immediately by explicit calculations and need not be derived from suitable Limiting Absorption Principles; for details, we refer to the earlier article~\cite{own_cubic}.

The study of conditions guaranteeing uniqueness of solutions of~\eqref{eq_linH} in $X_1$ involves the characterization of Herglotz waves in $X_1$, which are all multiples of $\tilde{\Psi}_\mu$. As in~\cite{own_cubic}, inspired by the analysis of the so-called far field of solutions of Helmholtz equations in scattering theory, we impose asymptotic conditions governing the leading-order contribution of $u(x)$ as $|x| \to \infty$. For $\tau \in (0, \pi)$, we introduce 
\begin{align*}
	\mathcal{R}_{\mu}^{\tau} [f] = \Psi_\mu \ast f + \cot(\tau) \: \tilde{\Psi}_\mu \ast f
	= \frac{\sin(|\,\cdot\,| \sqrt{\mu} + \tau)}{4\pi \sin(\tau) \: |\,\cdot\,|} \ast f.
\end{align*}
Then, using the above Lemma~\ref{lem_convol}, one obtains:
\begin{lem}[\cite{own_cubic}, Corollary~5]\label{lem_linH}
Let $\tau \in (0, \pi)$ and $\mu > 0$. Then the operator $\mathcal{R}_{\mu}^{\tau}: X_3 \to X_1$ is well-defined, linear and compact. Moreover, given $f \in X_3$, we have $u = \mathcal{R}_{\mu}^{\tau} [f]$ if and only if $u \in C^2_{\mathrm{loc}}$ with
\begin{align*}
	- \Delta u - \mu u = f \quad \text{on } \R^3,
	\qquad
	u(x) = c \cdot \frac{\sin(|x| \sqrt{\mu} + \tau)}{|x|} + O\left( \frac{1}{|x|^2} \right)
	\quad \text{as } |x| \to \infty
\end{align*} 
for some $c \in \R$, and in this case $c = \frac{1}{\sin(\tau)} \: \sqrt{\frac{\pi}{2}} \: \hat{f}(\sqrt{\mu})$.
\end{lem}
Handling the case of far field conditions with $\tau = 0$ is somewhat more delicate since the existence of the solution $\tilde{\Psi}_\mu$ (which satisfies exactly this condition) excludes an analogous uniqueness statement. For proving Theorem~\ref{thm_poly}, the following setting is suitable. First, by the Hahn-Banach Theorem, we define continuous linear functionals $\alpha^{(\mu)}, \beta^{(\mu)} \in X'_1$ with the property that, for $u \in X_1$ with
\begin{align*}
	u(x) = \alpha_u \cdot \tilde{\Psi}_\mu(x) + \beta_u \cdot \Psi_\mu(x) + O\left( \frac{1}{|x|^2} \right)
	\quad \text{as } |x| \to \infty,
\end{align*}
we have $\alpha^{(\mu)}(u) = \alpha_u$ and $\beta^{(\mu)}(u) = \beta_u$, cf. \cite{own_cubic}, equation~(13) and the following explanations. Then, the following analogue of Lemma~\ref{lem_linH} holds.

\begin{lem}\label{lem_linH-0}
Given $f \in X_3$, we have $u = \mathcal{R}_{\mu}^{\pi/2} [f] + (\alpha^{(\mu)}(u) + \beta^{(\mu)}(u)) \cdot \tilde{\Psi}_\mu$ if and only if $u \in C^2_{\mathrm{loc}}$ with
\begin{align*}
	- \Delta u - \mu u = f \quad \text{on } \R^3,
	\qquad
	u(x) = c \cdot \frac{\sin(|x| \sqrt{\mu})}{|x|} + O\left( \frac{1}{|x|^2} \right)
	\quad \text{as } |x| \to \infty
\end{align*} 
for some $c \in \R$. In this case, $\beta^{(\mu)}(u) = 0$.
\end{lem}
These results will allow to handle the nonlinear Helmholtz equations in~\eqref{eq_stationary_2}; for the proofs, we refer to the corresponding parts of~\cite{own_cubic}. Nonlinear Schrödinger equations such as
\begin{align}\label{eq_linS}
	- \Delta u + \mu u = f			\qquad \text{on } \R^3
\end{align}
for some $\mu > 0$ can also be discussed in a similar setting, which is certainly neither optimal nor most elegant but perfectly suitable for our purpose as another analogue of Lemma~\ref{lem_convol}. 
\begin{lem}\label{lem_linS}
Let $\mu > 0$. Then the operator
\begin{align*}
	\mathcal{P}_{\mu}: X_3 \to X_1, \quad
	f \mapsto \frac{\mathrm{e}^{-|\,\cdot\,| \sqrt{\mu}}}{4 \pi | \, \cdot \, |} \ast f
\end{align*}
is well-defined, linear and compact. Moreover, given $f \in X_3$, we have $u := \mathcal{P}_\mu [f]	\in X_3 \cap C^2_\mathrm{loc}(\R^3)$, and $u$ is a solution in $X_1$ of 
	\begin{align*}
			- \Delta u + \mu u = f  \quad \qquad \text{on } \R^3.
	\end{align*}		
\end{lem}
For details on the proof, which is similar to that of Lemma~\ref{lem_convol} but with less difficulties due to the strongly localized kernel, cf.~\cite{myDiss},~Lemma~4.10. 

Let us remark that, in the Schrödinger case, we do not obtain a family of possible ``resolvent-type'' operators as $\mathcal{R}_1^\tau = \mathcal{R}_1 + \cot(\tau) \tilde{\mathcal{R}}_1$, $0 < \tau < \pi$, in the Helmholtz case. This is due to the fact that the homogeneous Schrödinger equation $- \Delta u + \mu u = 0$ has no smooth and localized nontrivial solution in $X_1$. In particular, a major consequence in our study of Klein-Gordon breathers is that we have to impose nondegeneracy of $w_0$ as an assumption rather than, as in the Helmholtz case, generate it by choosing an appropriate resolvent  $\mathcal{R}_1^\tau$, as will be done in~\eqref{eq_assumptions},~\eqref{eq_nondegenerate,1} below.

\section*{Acknowledgements}

Funded by the Deutsche Forschungsgemeinschaft (DFG, German Research Foundation) – Project-ID 258734477 – SFB 1173. 

The construction of breather solutions of a similar wave-type equation is a major result of the author's dissertation thesis and can be found partly verbatim in~\cite[Chapter 4]{myDiss}. 
Special thanks goes in particular to my PhD advisor Dr. Rainer Mandel who encouraged me to work on this topic and provided advice whenever asked.

\bibliographystyle{abbrv}	
\bibliography{Literatur} 
  
\end{document}